\documentclass[12pt]{article}
\usepackage{latexsym}
\usepackage{amssymb}
\usepackage{amsmath}
\newcommand{\eproof}{\mbox{\ }\hfill $\Box$ \par \vskip 10pt}

\newtheorem{Theorem}{Theorem}
\newtheorem{lemma}{Lemma}
\newtheorem{prop}{Proposition}
\newtheorem{rem}{Remark}

\begin{document}

\title{Periodic Nonlinear Schr\"odinger Equation with Application to Photonic Crystals}

\author{
{\sc  Alexander Pankov}\\
Mathematics Department\\
College of William and Mary\\
Williamsburg, VA 23187--8795\\
e-mail: {\tt pankov@member.ams.org} }

\date{}

\maketitle

\begin{abstract} We present basic results, known and new, on nontrivial  solutions of periodic stationary nonlinear
Schr\"odinger equations. We also sketch an application to
nonlinear optics and discuss some open problems.
\end{abstract}

\setcounter{section}{-1}
\section{Introduction}

We consider the nonlinear Schr\"odinger equation (NLS)
\begin{equation}\label{e1}
-\Delta u +V(x)u=\pm f(x,u), \quad x\in \mathbb{R}^N,
\end{equation}
with superlinear subcritical nonlinearity $f$, $f(x,0)=0$. Our basic assumption is periodicity of
$V(x)$ and $f(x,u)$ in $x$ with respect to some lattice $\Lambda\subset \mathbb{R}^N$. For the sake
of simplicity in the following we assume that $\Lambda = \mathbb{Z}^N$ is the standard integer lattice in
$\mathbb{R}^N$, i.e. the functions $V(x)$ and $f(x,u)$ are 1-periodic in each variable $x_i$,
$x=(x_1, x_2, \ldots , x_N)$. We are interesting in nontrivial solutions of equation (\ref{e1}) such that
$u(\infty)=0$.

Problems of such kind appear in many applications. For instance, in condensed matter physics, the states
of atomic gas are described by wave functions $\psi(t,x)$ that satisfy the Gross-Pitaevski
equation (see, e.g., \cite{bec})
\[
i\psi_t=-\Delta\psi +V_0(x)\psi-\alpha |\psi|^2\psi
\]
(the so-called mean field approximation). Condensed states confined in a finite region
(Bose-Einstein condensates) correspond to the wave
functions $\psi=u(x)\cdot\exp(-i\omega t)$ where the amplitude $u(x)$ vanishes at infinity, i.e. standing
waves. Here
$V_0(x)$ is the external electric potential which can be periodic. Making use the standing wave Ansatz, we arrive
at the equation of the form (\ref{e1}). Another field of applications, gap solitons in photonic crystals,
will be discussed in Section~\ref{gs} below.

There is a number of papers dealing with periodic NLS and its solutions vanishing at infinity. First, the
case when the operator $-\Delta +V$ is positive definite was considered. In \cite{pank:umj} the author
has proved an existence theorem using the Nehari variational principle (see, e.g., \cite{will})
and concentrated compactness results \cite{lions}. (Even more general asymptotically periodic case was treated
in that paper). Later  P.~H.~Rabinowitz \cite{rab:per} has obtained the existence of nontrivial solutions under
less restrictive assumptions on the nonlinearity $f$. Moreover, in \cite{cz-rab} V.~Coti Zelati and P.~H.~Rabinowitz
have proved the existence of infinitely many multi-bump solutions, imposing some condition of general position.

In general, it is known (see, e.g.,\cite{r-s}) that the spectrum of the self-adjoint operator $-\Delta+V$
in ${\mathrm L}^2(\mathbb{R}^N)$ is purely continuous and may contain gaps, i.e. open intervals free of spectrum.
So, it is possible that 0 is in a finite spectral gap and the operator $-\Delta+V$ is not positive definite.
In this case first existence results (under very strong assumptions on the nonlinearity) were found in
\cite{ala-li1, je}. Later C.~Troestler and M.~Willem \cite{t-will} and W.~Kryszewski and A.~Szulkin \cite{k-sz}
have proved the existence of nontrivial solutions under much more natural conditions. Their proofs are based on
a generalized linking theorems applied to the corresponding action functional. The original proofs of such linking
results required new degree theories in order to overcome the lack of compactness and find at least one Palais-Smale
sequence. Moreover, in \cite{k-sz} it is shown that there exists infinitely many geometrically distinct solutions,
provided the nonlinearity is odd. It is not known whether these solutions are multibumps.
In \cite{ba-ding},  T.~Bartsch and Y.~Ding considered the case when 0 is a boundary point of the spectrum
and, in addition, gave a simpler proof the generalized linking theorem. In \cite{pank:nls}, the author and K.~Pfl\"uger
suggested another approach to the periodic NLS based on periodic approximations
and standard critical point theorems. This approach goes back to
P.~H.~Rabinowitz \cite{rab:per} who considered the positive definite case. Later the same idea was applied to other
problems \cite{pank:kp1, pank:kp2, pank:lat}. See also \cite{ack, ala-li1, ala-li2, je, bu-je-st, li-sz, wi-z}
for some other results on NLS.

Now we list the basic assumptions. Let
$$F(x,u)=\int_0^u f(x,s)ds.$$
\begin{description}
\item[$(i)$] {\em $V\in ^{\infty}(\mathbb{R}^N)$  is a 1-periodic function in each $x_i$,
              $i=1,\dots , N$.}

\item[$(ii)$] {\em $f:\mathbb{R}^N\times\mathbb{R} \to \mathbb{R}$ is a Carath\'eodory function\footnote{This means that
              $f(x,u)$ is measurable in $x\in \mathbb{R}^N$ for all $u\in \mathbb{R}$ and
               continuous in $u\in \mathbb{R}$ for a.e. $x\in \mathbb{R}^N$.} that is 1-periodic in each $x_i$,
              $i=1,\dots , N$.}

\item[$(iii)$] {\em $|f(x,u)|\leq C(1+|u|^{p-1})$ where $2<p<2^*$,  $2^*:= 2N/(N-2)$ if $N\geq 3$
       and $2^*:=+\infty$ if $N=1, 2$.}

\item[$(iv)$] $${\mathrm{ess}}\sup_{x\in \mathbb{R}^N} |f(x,u)|=o(|u|), \quad u\to 0.$$

\item[$(v)$] {\em There exists $q >2$ such that}
    $$0< q F(x,u) \leq uf(x,u), \quad u\neq 0.$$
\end{description}

Assumption $(i)$ implies that the operator $-\Delta+V$ with the domain ${\mathrm C}_0^{\infty}(\mathbb{R}^N)$
is essentially self-adjoint in $L^2(\mathbb{R}^N)$. Denote by $L$ its self-adjoint extension and
$\sigma (L)$ the spectrum of $L$. The last assumption is the following
\begin{description}
\item[$(vi)$] {\em $0$ is not in $\sigma (L)$. In addition, in the case of `$-$' in equation (\ref{e1})
             we assume that $\inf\sigma (L) < 0$.}
\end{description}
This assumption implies that there is a maximal open interval $(-\alpha_-, \alpha_+)$, with $\alpha_{\pm}>0$,
free of the spectrum
(gap in the spectrum).
In the case of `$+$' sign it is
possible that $\alpha_-=\infty$. However, in the case of `$-$' sign both $\alpha_-$ and $\alpha_+$
are finite, i.e. 0 is in a
finite gap. We set $\alpha:=\min[\alpha_+,\alpha_-]$. Note that $\alpha$ is the distance from 0 to $\sigma (L)$.

The power nonlinearity $f(x,u)=h(x)|u|^{p-2}u$, with positive 1-periodic $h\in  L^{\infty}(\mathbb{R}^N)$
and $2<p<2^*$, satisfies all the assumptions. The cubic nonlinearity ($p=4$) is admissible if $N=1, 2$ or $3$. The
last nonlinearity is important for many applications (see, e.g., Section~\ref{gs}). Note that in the case $N=1$
assumption $(iii)$ can be replaced by the following one. For any $h>0$ there exists a constant $C(h)>0$ such that
$${\mathrm{ess}}\sup_{x\in \mathbb{R}^N} |f(x,u)|\leq C(h) \quad \text{if\ }|u|\leq h.$$

Assumptions $(iii)$ and $(v)$ mean that we consider subcritical superlinear nonlinearities. For some results for  the critical periodic NLS
we refer to \cite{chab}. Another important problem concerns the existence of solutions to (\ref{e1}) vanishing at infinity in the case of
asymptotically linear nonlinearity (see a remark in Section~\ref{gs}). Some results in this direction can be found in \cite{li-sz}, but the problem
is still not well understood.

As it follows from $(v)$, the nonlinearity $f(x,u)$ does not change sign when $x$ changes. In the case of periodic NLS with sign changing nonlinearity
the existence of vanishing at infinity solutions is a completely open problem that may have interesting consequences in the theory of gap solitons
(see discussion in Section~\ref{gs}).

Now we summarize the principal result.

\begin{Theorem}\label{t0}
Under assumptions $(i)$--$(vi)$ equation (\ref{e1}) has a nontrivial weak solution
$u\in \mathrm{H}^1(\mathbb{R}^N)$. This solution is continuous and decays exponentially at infinity.
\end{Theorem}

In  subsequent sections we present the proof and discuss some
related results. We give also an application to the problem of
existence of gap solitons in photonic crystals.

Throughout the paper we deal basically with the case of `$+$' sign in (\ref{e1}) and only point out minor
changes needed to cover the case of `$-$' sign.

\vspace{2ex}

{\bf Acknowledgments.} This article was prepared during author's
staying at Texas A\&M University and College of William \& Mary as
a visiting professor. The author also thanks T.~Bartsch,
P.~Kuchment, P.~H.~Rabinowitz, A.~Szulkin and I.~Spitkovsky for many
interesting discussions. The work is partially supported by NATO, grant 970179.

\section{Abstract critical point theorems}\label{s1}

Let $H$ be a real  Hilbert space with the inner product $(\cdot ,\cdot )$ and the norm $\|\cdot\|$,
and $J$ a functional on $H$ of
the class $C^1(H)$. Denote by $J'(u)$ the Gateaux derivative (gradient) of $J$ at the point $u\in H$.
Recall that $u\in H$ is a {\em critical point} of $J$ if $J'(u)=0$ and $c=J(u)$ is the corresponding
{\em critical value}.
We shall use the following Palais-Smale condition:
\begin{description}
\item[$(PS)$] {\em Any sequence $u_n\in H$ such that $J(u_n)$ is bounded and $J'(u_n)\to 0$, i.e.
{\em Palais-Smale sequence\/}, has a convergent
               subsequence.}
\end{description}




Now we present the Linking Theorem (see \cite{rab:cbms, will} for the proof).
Let $H=Y\oplus Z$, $\rho > r >0$ and let $z\in Z$ be such that
$\|z\| = r$. Define
$$M:=\{u=y+\lambda z: \|u\|\leq\rho, \lambda \geq 0, y\in Y\},$$
$$M_0:=\partial M=\{u=y+\lambda z : y\in Y,\|u\|=\rho \mbox{\ and\ }\lambda\geq 0,
\mbox{\ or\ } \|u\|\leq\rho \mbox{\ and\ } \lambda =0\},$$
$$N:=\{u\in Z : \|u\|=r\}.$$
We say that the functional $J$ on $H$ possesses the {\em Linking Geometry} if
$$b:=\inf_{u\in N} J(u) > a:=\sup_{u\in M_0} J(u).$$

\begin{Theorem}\label{link}

Suppose $J$ satisfies $(PS)$ and possesses the Linking Geometry with $\mbox{\rm dim}\,Y<\infty$.
Let
$$c:=\inf_{\gamma\in\Gamma}\sup_{u\in M} J(\gamma(u))$$
where
$$\Gamma:=\{\gamma\in C(M, H) : \gamma =\mathrm{id} \mbox{\ on\ } M_0\}.$$
Then $c$ is a critical value of $J$ and $b\leq c\leq\sup_{u\in M}J(u)$.
\end{Theorem}

Note that in the case $Y=\{0\}$ Theorem~\ref{link} reduces to the Mountain Pass Theorem (see, e.g.,
\cite{rab:cbms, will}).

Sometimes one can drop the assumption $\mbox{dim}\,Y<\infty$ in the Linking Theorem. The following result was
obtained in \cite{b-rab}.

\begin{Theorem}\label{link1}
Assume that the functional $J$ possesses the Linking Geometry  and satisfies
the following assumptions
\begin{description}
\item[$(j)$] $J(u)=\frac{1}{2}(Au,u)+b(u)$ where $u=u_1+u_2\in Y\oplus Z$, $Au=A_1u_1+A_2u_2$ $A_1:Y\to Y$
             and $A_2:Z\to Z$ are linear bounded self-adjoint operators;
\item[$(jj)$] $b$ is weakly continuous and uniformly differentiable on bounded sets;
\item[$(jjj)$] Any sequence $u_n\in H$ such that $J(u_n)$ is bounded from above and $J'(u_n)\to 0$ is bounded.
\end{description}
Then $c$ defined in Theorem~\ref{link} is a critical value of $J$ such that $b\leq c\leq\sup_{u\in M}J(u)$.
\end{Theorem}

Note that in the case of Theorem~\ref{link1} the functional $J$ satisfies the Palais-Smale condition.

\section{Variational setting and periodic problem}\label{s2}

Throughout this section we will suppose  assumptions $(i)$--$(vi)$
to be valid.

Let $Q_k$ be the cube in $\mathbb{R}^N$, with the edge length $k$, centered at the origin. Denote by $L_k$ a unique
self-adjoint operator in $L^2(Q_k)$ generated by $-\Delta +V$ with periodic boundary conditions. We can (and will
consider) this space as the space of all $k$-periodic functions from $L^2_{\mathrm{loc}}(\mathbb{R}^N)$. Note that
$\sigma(L_k)$ is discrete and is a subset of $\sigma(L)$ (see, e.g., \cite{r-s}). Let $E_k=H^1_{\mathrm{per}}(Q_k)$ be
the subspace of $H^1_{\mathrm{loc}}(\mathbb{R}^N)$ that consists of $k$-periodic functions. This is just the form
domain of $L_k$. Endowed with the standard $H^1(Q_k)$ norm, $E_k$ is a Hilbert space. We also set
$E=H^1(\mathbb{R}^N)$. This space coincides with the form domain of $L$.

Let $E_k^+$ (respectively, $E_k^-$) be the intersection of $E_k$ with the positive (respectively, negative)
spectral subspace of $L_k$. Similarly we introduce the subspaces $E^+$ and $E^-$ in $E$. It is easily seen that
$E_k^+=(E_k^-)^{\perp}$ and $E^+=(E^-)^{\perp}$. Thus, $E_k= E_k^-\oplus E_k^+$ and $E_k= E_k^-\oplus E_k^+$.
These decompositions define the orthogonal projections $P_k^{\pm}$ onto $E_k^{\pm}$ and
$P^{\pm}$ onto $E^{\pm}$. For the sake of convenience, we introduce new inner products $(\cdot ,\cdot)_k$ and
$(\cdot ,\cdot)$ in $E_k$ and $E$, with corresponding norms $\|\cdot\|_k$ and $\|\cdot\|$, so that
$$
\int_{Q_k}(|\nabla u|^2+V(x)|u|^2)\,dx=-\|u\|_k^2 \quad\mbox{for}\quad u\in E^-_k\,,
$$
$$
\int_{Q_k}(|\nabla u|^2+V(x)|u|^2)\,dx=\|u\|_k^2 \quad\mbox{for}\quad u\in E^+_k\,
$$
and
$$
\int_{\mathbb{R}^N}(|\nabla u|^2+V(x)|u|^2)\,dx=-\|u\|^2 \quad\mbox{for}\quad u\in E^-\,,
$$
$$
\int_{\mathbb{R}^N}(|\nabla u|^2+V(x)|u|^2)\,dx=\|u\|^2 \quad\mbox{for}\quad u\in E^+\,.
$$
These norms are equivalent to the standard $H^1$-norms independently of $k$, i.e. there
exists a constant $C>0$, independent of $k$, such that
$$
\alpha^{1/2}\|u\|_{H^1(Q_k)}\leq \|u\|_k\leq C \|u\|_{H^1(Q_k)}
$$
and
$$
\alpha^{1/2}\|u\|_{H^1(\mathbb{R}^N)}\leq \|u\|\leq C \|u\|_{H^1(\mathbb{R}^N)}.
$$
On the subspaces $E_k^{\pm}$ and $E^{\pm}$ the constant  $\alpha$ can be replaced by $\alpha_{\pm}$.
Moreover, the constant $C$ can be chosen independent of $V$ when $V$ runs a bounded subset of
$L^{\infty}(\mathbb{R}^N)$.

Consider the action functionals
\begin{align*}
J_k(u)&=\frac{1}{2}(\|P_k^+u\|^2_k-\|P^-_ku\|^2)\mp\int_{Q_k}F(x,u)\, dx=\\
      &=\frac{1}{2}\int_{Q_k}(|\nabla u|^2+V(x)|u|^2 \mp F(x,u))\, dx\\
\intertext{and}
J(u)&=\frac{1}{2}(\|P^+u\|^2-\|P^-u\|^2)\mp\int_{\mathbb{R}^N}F(x,u)\, dx=\\
      &=\frac{1}{2}\int_{\mathbb{R}^N}(|\nabla u|^2+V(x)|u|^2 \mp F(x,u))\, dx
\end{align*}
on $E_k$ and $E$, respectively. From the standard Sobolev embedding theorem (see, e.g., \cite{adams, maz})
it follows that these functionals are well defined. One can also check that they are of the class $C^1$
and  the derivatives are given by the formulas
\begin{align*}
\langle J'_k(u),v\rangle&=(P_k^+u,v)_k-(P^-_ku,v)_k\mp\int_{Q_k}f(x,u)v\, dx=\\
      &=\int_{Q_k}(\nabla u\cdot\nabla v+V(x)uv \mp f(x,u)v)\, dx, \quad \forall v\in E_k\\
\intertext{and}
\langle J'(u),v\rangle&=(P^+u,v)-(P^-u,v)\mp\int_{\mathbb{R}^N}f(x,u)v\, dx=\\
      &=\int_{\mathbb{R}^N}(\nabla u\cdot\nabla v+V(x)uv \mp f(x,u)v)\, dx \quad \forall v\in E.
\end{align*}
Therefore, critical points of these functionals are weak solutions of (\ref{e1}), $k$-periodic in the case of
$J_k$ and vanishing, in a sense, at infinity in the case of $J$.

Recall that the gradient $\nabla J_k(u)\in E_k$ (resp., $\nabla J(u)\in E$) is defined by $(\nabla J_k(u), v)_k=\langle J'_k(u),v\rangle$
(resp., $(\nabla J(u), v)_k=\langle J'(u),v\rangle$) for all $v\in E_k$ (resp., $v\in E$).

\begin{lemma}\label{l2.1} There exists a constant $C>0$  independent on $k$
such that for any nontrivial critical points $u_k\in E_k$ of $J_k$ and $u\in E$ of $J$, with critical
values $c_k$ and $c$ respectively, we have
$$\|u_k\|_k\leq C(|c_k|^{1/2}+|c_k|^{1/p'}),\quad \|u\|\leq C(|c|^{1/2}+|c|^{1/p'})$$
where $p':=p/(p-1)$ is the conjugate exponent.
\end{lemma}

{\em Proof\/}. From $c_k=J_k(u_k)$ and $J'_k(u_k)=0$ we obtain
\begin{align}\label{e2.1}
|c_k|&= |J_k(u_k)-\frac{1}{2}\langle J'_k(u_k),u_k\rangle| \geq\notag\\
      &\geq\left(\frac{1}{2}-\frac{1}{q}\right)\int_{Q_k}f(x,u_k)u_k\, dx.
\end{align}
From assumptions $(iii)$--$(v)$ we have the inequalities
\begin{equation}\label{e2.2}
|f(x,u)|^2\leq Cuf(x,u)\quad \mbox{if}\quad |u|\leq 1,
\end{equation}
\begin{equation}\label{e2.3}
|f(x,u)|^{p'}\leq C|u|^{(p-1)(p'-1)}f(x,u)=Cuf(x,u)\quad \mbox{if}\quad |u|\geq 1.
\end{equation}
Let
$$B_k=\{x\in Q_k\,|\,|u_k(x)\leq 1\}$$
Then from (\ref{e2.1})--(\ref{e2.3}) we obtain
\begin{equation}\label{e2.4}
|c_k|\geq C'\left(\int_{B_k}|f(x,u_k(x))|^2\, dx +\int_{Q_k\setminus B_k}|f(x,u_k(x))|^{p'}\, dx\right)
\end{equation}
where $C'=(q-2)/(2qC)$. Hence
$$I_1:=\int_{B_k}|f(x,u_k(x)|^2\, dx\leq |c_k|/C',$$
$$I_2:=\int_{Q_k\setminus B_k}|f(x,u_k(x))|^{p'}\, dx\leq |c_k|/C'.$$
Let $y_k=P^-_ku_k$ and  $z_k=P^+_ku_k$. From $\langle J'_k(u_k),y_k\rangle =0$ and H\"older's inequality
we see that
$$\|y_k\|^2_k=-\int_{Q_k}f(x,u_k)y_k\, dx\leq I^{1/2}_1|y_k|_2+I^{1/2}_2|y_k|_{p'}$$
where $|\cdot |_r$ stands for the standard $L^r$-norm. Then the Sobolev embedding theorem
implies
$$\|y_k\|^2_k\leq C''\left(I_1^{1/2}+I_2^{1/p'}\right)\|y_k\|_k.$$
Clearly, the same argument works for $z_k$ and we obtain
$$\|z_k\|^2_k\leq C''\left(I_1^{1/2}+I_2^{1/p'}\right)\|z_k\|_k.$$
These two inequalities imply the first estimate of the lemma. Similarly, we obtain the second one.
\eproof

\begin{rem}\label{r2.a}
{\em Similarly, we obtain the following more precise estimates}
$$
\alpha^{1/2}\|u_k\|_{H^1(Q_k)}\leq C(|c_k|^{1/2}+|c_k|^{1/p'}),
\quad \alpha^{1/2}\|u\|_{H^1(\mathbb{R}^N)}\leq C(|c|^{1/2}+|c|^{1/p'})\,.
$$
\end{rem}

\begin{lemma}\label{l2.2}
There exist constants $\varepsilon_1 >0$ and $\varepsilon_2 >0$ independent on $k$ such that
for any nontrivial critical points $u_k\in E_k$ of $J_k$ and $u\in E$ of $J$ we have
$\|u_k\|_k \geq \varepsilon_1$, $\|u\| \geq \varepsilon_1$, $\mp J_k(u_k)\geq\varepsilon_2$
and $\pm J(u)\geq\varepsilon_2$.
\end{lemma}

{Proof\/}. Consider `$+$' sign case. Assumptions $(iii)$ and $(iv)$ implies immediately that for any
$\varepsilon >0$ there exists $C_{\varepsilon}>0$ such that
\[|f(x,t)|\leq \varepsilon |t|+C_{\varepsilon}|t|^{p-1}.\]
Since $J'_k(u_k)=0$, then $\langle J'_k(u_k),v_k\rangle =0$ for any $v_k\in E_k$.
Take $v_k=u_k^+=P_k^+u_k$. Then we have
\begin{align*}
\|u_k^+\|_k^2&=\int_{Q_k} f(x,u_k)u_k^+\, dx \leq \\
         &\leq C(\varepsilon \|u_k\|_{L^2(Q_k)}\|u_k^+\|_{L^2(Q_k)}+
                       C_{\varepsilon}\|u_k\|^{p-1}_{L^p(Q_k)}\|u_k^+\|_{L^p(Q_k)})\\
         &\leq C(\varepsilon \|u_k\|_k\|u_k^+\|_k+
                       C_{\varepsilon}\|u_k\|_k^{p-1}\|u_k^+\|_k)
\end{align*}
where we have used the Sobolev embedding theorem. Similarly, taking $v_k=u_k^-=P_k^-u_k$ we obtain
\begin{align*}
\|u_k^-\|^2_k&=\int_{Q_k} f(x,u_k)u_k^-\, dx \leq \\
         &\leq C(\varepsilon \|u_k\|_k\|u_k^-\|_k+
                       C_{\varepsilon}\|u_k\|_k^{p-1}\|u_k^-\|_k).
\end{align*}
These two inequalities imply immediately
\[
\|u_k\|_k^2=\|u_k^+\|_k^2+\|u_k^-\|_k^2
         \leq C (\varepsilon \|u_k\|_k^2+
                       C_{\varepsilon}\|u_k\|_k^{p})\,.
\]
If $\varepsilon$ is small enough, we obtain
\[
\|u_k\|_k\geq \left(\frac{1-\varepsilon C}{C_{\varepsilon}}\right)^{1/(p-2)}=:\varepsilon_1 >0.
\]
The same argument, with the same choice of constants, works for $J$.
The estimates for critical values follow from Lemma~\ref{l2.1}.

The case of `$-$' sign is similar. \eproof

\begin{Theorem}\label{t2.1}The functional $J_k$ has a nontrivial critical point $u_k\in E_k$ with
critical value $c_k$. Moreover, there exist a constant $C>0$ independent of $k$ such that
$\|u_k\|_k\leq C$ and $0< \pm c_k\leq C$.
\end{Theorem}

{\em Proof\/}. We start with the Palais-Smale condition. Refining the proof of Lemma~\ref{l2.1}, one can show
that any Palais-Smale sequence for $J_k$ is bounded in $E_k$. Now, using the standard argument
\cite{rab:cbms, will} based on the compactness of Sobolev embedding, we conclude that such a sequence is compact
in $E_k$.

Consider `$+$' sign case. We show that the functional $J_k$ possesses the linking geometry with
$Y=E^-_k$ and $Z=E^+_k$.

For $z\in Z$ we have
$$
J_k(z)=\frac{1}{2}\|z\|^2_k-\int_{Q_k}F(x,z)\, dx.
$$
Assumptions $(iii)$ and $(iv)$ imply that for any $\varepsilon >0$ there exists a constant
$C_{\varepsilon} >0$ such that $0\leq F(x,u)\leq\varepsilon |u|^2+C_{\varepsilon}|u|^p$.
By the Sobolev embedding theorem we obtain
$$
\int_{Q_k}F(x,z)\, dx \leq C(\varepsilon \|z\|^2_k+C_{\varepsilon}\|z\|_k^p).
$$
Choosing an appropriate $\varepsilon$  we see that, for some $\delta >0$, $J_k(z)\geq\delta$ if
$\|z\|_k=r$ is small enough. Thus, we have chosen the sphere $N=\{u\in Z\,:\, \|u\|_k=r\}$.

Let us now find $M$. Fix any $z^0\in Z=E^+_k$, $\|z^0\|_k=1$. For $y+tz^0\in M$, $y\in Y=E^-_k$
we have
$$
J_k(y+tz^0)=\frac{1}{2}t^2-\frac{1}{2}\|y\|^2_k-\int_{Q_k}F(x,y+tz^0)\, dx.
$$
Assumptions $(iv)$ and $(v)$ imply that for any $\varepsilon >0$ there exists a constant
$C_{\varepsilon}>0$ such that $F(x,u)\geq -\varepsilon |u|^2 +C_{\varepsilon}|u|^q$. Hence,
$$
\int_{Q_k}F(x,y+tz^0)\, dx \geq -\varepsilon\|y\|^2_{L^2(Q_k)}-
\varepsilon\|tz^0\|^2_{L^2(Q_k)} +C_{\varepsilon}\|y+tz^0\|^q_{L^q(Q_k)}.
$$
Applying the Sobolev embedding theorem, we obtain
\begin{equation}\label{e2.5}
J_k(y+tz^0)\leq -\frac{1}{2}\|y\|^2_k+\varepsilon C\|y\|^2_k+\frac{1}{2}t^2 +
\varepsilon Ct^2 - C_{\varepsilon}\|y+tz^0\|^q_{L^q(Q_k)}.
\end{equation}
Consider the space $X=Y\oplus \mathbb{R}z^0$ equipped with the $L^q(Q_k)$-norm. Then
$y+tz^0\mapsto tz^0$ is a bounded projector in $X$. Since its norm is not less then 1, we see that
$$
\|y+tz^0\|_{L^q(Q_k)}\geq\|tz^0\|_{L^q(Q_k)}.
$$
Choosing $\varepsilon$ so small that $\varepsilon C=1/4$, we obtain from (\ref{e2.5})
$$
J_k(y+tz^0)\leq-\frac{1}{4}\|y\|^2_k +\frac{3}{4}t^2-C_{\varepsilon}C|t|^p
$$
(here we have used again the Sobolev embedding). Note that all the constants here are independent on $k$.
The last inequality implies that we can choose $\rho >0$ large enough such that $J_k\leq 0$ on $M_0$.
Moreover, $\sup_{M}J_k(u)\leq K:=\max(3t^2/4-C_{\varepsilon}C|t|^q)$ and $K$ is independent of $k$.
Applying Theorem~\ref{link}, we obtain the result.

The case `$-$' sign we need only to apply Theorem~\ref{link1} to the functional
$-J_k$. \eproof

\begin{rem}\label{r2.*}{\em If we drop assumption $(vi)$, i.e. allow 0 to be in the spectrum $\sigma(L)$, we still can use linking
argument to prove the existence of $k$-periodic solutions. This can be done as in the case of the Dirichlet problem in a bounded domain (see
\cite{will}, Section~2.4). However, in this case we have no uniform (with respect to $k$) bounds for the solution. }
\end{rem}

\begin{rem}\label{r2.b}{\em In the case `$-$', if we drop the requirement $\inf\sigma (L) < 0$ then there is no
nontrivial solution in $E_k$, as well as in $E$. This follows from the obvious fact that in this case the
functionals $-J_k$ and $-J$ have strict global maximum at the origin.}
\end{rem}

\section{Passage to the limit and finite action solutions}\label{s3}

Now we are going to pass to the limit as $k\to\infty$ to obtain a finite action solution $u\in E$.
The principal difficulty is to show that the limit is a non-zero function. With this aim we need
some concentration compactness argument.
\begin{lemma}\label{l3.1}
Let $Q_{(n)}$ be the cube with the edge length $l_n\to\infty$ centered at the origin and $K_r(\xi)$ the cube
centered at $\xi$ with the edge length $r$. Let $u_{(n)}\in H^1_{\mathrm{loc}}(\mathbb{R}^N)$ be
a sequence of $l_n$-periodic functions such that $\|u_{(n)}\|_{H^1(Q_{(n)})}\leq C$. Assume that
there exists $r>0$ such that
$$
\lim_{n\to\infty}\left(\sup_{\xi}\int_{K_r(\xi)}|u_{(n)}|^2\, dx\right)=0.
$$
Then $\|u_{(n)}\|_{L^s(Q_{(n)})}\to 0$ for all $s\in (2, 2^*)$.
\end{lemma}

{\em Proof\/}. Denote by $Q'_{(n)}$ and $Q''_{(n)}$ the cubes with the sizes $l_n+r$ and $l_n+2r$,
respectively, centered at 0. Choose a cut-off function $\chi_n\in C_0^{\infty}(\mathbb{R}^N)$
such that $\mbox{supp}\,\chi_n\subset Q''_{(n)}$, $0\leq \chi_n(x)\leq 1$, $\chi_n = 1$ on $Q'_{(n)}$
and $|\nabla\chi_n(x)\|\leq C$ where $C>0$ does not depend on $n$. Evidently, such a function exists.

Let us set $v_n=\chi_n u_{(n)}$. For $n$ large enough, it is easy to verify, using the periodicity of $u_{(n)}$,
that
$$
\sup_{\xi}\int_{K_r(\xi)}|v_n|^2\, dx = \sup_{\xi}\int_{K_r(\xi)}|u_{(n)}|^2\, dx
$$
and the sequence $v_n$ is bounded in $H^1(\mathbb{R}^N)$. By the well-known lemma of P.~L.~Lions
(see \cite{lions}, Lemma~I.1, and \cite{will}, Lemma~1.21), $v_n\to 0$ in $L^s(\mathbb{R}^N)$,
$s\in (2, 2^*)$. Since $v_n=u_{(n)}$ on $Q_{(n)}$, we obtain the required. \eproof

\begin{lemma}\label{l3.2}
Let $u_k\in E_k$ be a sequence such that $\|u_k\|_k$ is bounded and $J'_k(u_k)\to 0$. Then, passing to a subsequence
still labelled by $k$, either
\begin{description}\item[$1^0$] $\|u_k\|_k\to 0$ as $k\to \infty$
\end{description}
or
\begin{description}
      \item[$2^0$] there exist a sequence of points $\xi_k\in\mathbb{R}^N$ and positive numbers $r$ and $\eta$
       such that
$$
\lim_{k\to\infty}\int_{K_r(\xi_k)}|u_k|^2\, dx \geq\eta.
$$
\end{description}
\end{lemma}

{\em Proof\/}. Assume that $2^0$ does not hold (along a subsequence still labelled by $k$). Since $J'_k(u_k)\to 0$,
there exists a sequence $\varepsilon_k\to 0$ such that $\langle J'_k(u_k), v\rangle\leq\varepsilon_k\|v\|_k$ for any
$v\in E_k$. Taking $v=P^+_ku_k$ and $v=P^-_ku_k$ we obtain, as in the proof of Lemma~\ref{l2.2},
$$
\|u_k\|^2_k\leq \varepsilon\|u_k\|^2_k +C_{\varepsilon}\|u_k\|^p_{L^p(Q_k)}+
C\varepsilon_k\|u_k\|_k\,.
$$
By Lemma~\ref{l3.1}, $\|u_k\|^p_{L^p(Q_k)}\to 0$ and we obtain the required.\eproof

\begin{Theorem}\label{t3.1}
Let $u_k\in E_k$ be a sequence of nontrivial critical points of $J_k$ such that the sequence
$c_k=J_k(u_k)$ is bounded. Then there exists a nontrivial solution $u\in E=H^1(\mathbb{R}^N)$
of equation (\ref{e1}). Moreover, there exist a sequence $b_k\in \mathbb{Z}^N$ such that  $u_k(x+b_k)$
converges to $u$  in $H^1_{\mathrm{loc}}(\mathbb{R}^N)$, up to the passage to a subsequence.
\end{Theorem}

{\em Proof\/}. By Lemma~\ref{l2.1}, the norms $\|u_k\|_k$ are bounded. Lemma~\ref{l2.2} implies that
$\|u_k\|-k\geq \varepsilon_1 >0$. Therefore, case $1^0$ of Lemma~\ref{l3.2} is not possible. From
Lemma~\ref{l3.2}, $2^0$, we obtain that, along a subsequence,
\[
\|u_k\|^2_{L^2(K_r(\xi_k))}\geq\eta/2.
\]
There exists a sequence of integer vectors $b_k$ such that the sequence of cubes $K_r(\xi_k-b_k)$ is
confined in a bounded region. Set $\tilde u_k(x)=u_k(x+b_k)$. We have
\begin{equation}\label{e3.1}
\|u_k\|^2_{L^2(Q_{k_0})}\geq\eta/2
\end{equation}
for  $k\geq k_0$ and some $k_0$.
Since $V$ and $f$ are $\mathbb{Z}^N$-periodic,
$\tilde u_k$ are critical points of $J_k$ and the norms $\|\tilde u_k\|_k$ are bounded. Passing to a subsequence,
one can assume that the sequence $\tilde u_k$ converges weakly in  $H^1_{\mathrm{loc}}(\mathbb{R}^N)$ to
a function $u$. Since $\|\tilde u_k\|_k$ is bounded,
$u\in H^1(\mathbb{R}^N)$. By the compactness of Sobolev embedding, $\tilde u_k\to u$ in
$L^{s}_{\mathrm{loc}}(\mathbb{R}^N)$, $2\leq s <2^*$, and, hence, $f(x,\tilde u_k)\to f(x,u)$
in $L^1_{\mathrm{loc}}(\mathbb{R}^N)$. Now for any test function $\varphi\in C_0^{\infty}$ we have
\begin{align*}
\langle J'(u),\varphi\rangle&=\int_{\mathbb{R}^N}\left(\nabla u\nabla\varphi +Vu\varphi - f(x,u)\varphi\right)\, dx = \\
         &=\lim_{k\to\infty}\int_{\mathbb{R}^N}\left(\nabla \tilde u_k\nabla\varphi +V\tilde u_k\varphi
                 - f(x,\tilde u_k)\varphi\right)\, dx =\\
         &= 0\,.
\end{align*}
Therefore, $u$ is a weak solution of (\ref{e1}) and $u\neq 0$ in view of (\ref{e3.1}).

Now let us prove the last statement. Let $\chi\in C_0^{\infty}(\mathbb{R}^N)$ and
$\mbox{supp}\,\chi$ is contained in
an open bounded set $\Omega$.
It is easy to verify the following
identity
$$
L(\chi \tilde u_k-\chi u)=h^1_k+h^2_k+h^3_k
$$
where
\begin{align*}
h^1_k&=\chi [f(x,\tilde u_k)-f(x,u)]\,, \\
h^2_k&=-\nabla\chi\nabla(\tilde u_k-u)\,,\\
h^3_k&=-(\Delta\chi)(\tilde u_k-u)\,.
\end{align*}
We see that $h^1_k\to 0$ in $L^{p-1}(\Omega)$, $h^2_k\to 0$ weakly in $L^2(\Omega)$ and
$h^3_k\to 0$ in $L^2(\Omega)$. Since the space $L^s(\Omega)$, $2\leq s <2^*$, is compactly
embedded into $H^{-1}(\mathbb{R}^N)$, then $h^1_k+h^2_k+h^3_k\to 0$ in the last space. Assumption
$(vi)$ implies that the operator $L$ considered as a bounded linear operator
from $H^1(\mathbb{R}^N)$ into $H^{-1}(\mathbb{R}^N)$ has a bounded inverse operator.
Therefore, $(\chi \tilde u_k-\chi u)\to 0$ in $H^1(\mathbb{R}^N)$ and we are done. \eproof

Combining Theorems~\ref{t2.1} and \ref{t3.1}, we obtain

\begin{Theorem}\label{t3.2}
Equation (\ref{e1}) possesses a nontrivial solution in $H^1(\mathbb{R}^N)$.
\end{Theorem}

\section{Ground states}\label{s4}

We say that a nontrivial solution $u\in E$ (resp., $u\in E_k$) of equation (\ref{e1}) is a {\em ground state}
(resp., a {\em periodic ground state}) if the corresponding critical value of the functional $\pm J$ (resp.,
$\pm J_k$) is minimal possible among all nontrivial critical values of that functional. Since the functional $\pm J_k$
satisfies the Palais-Smale condition, the existence of periodic ground states follows immediately from Theorem~\ref{t2.1}
and Lemma~\ref{l2.2}. The case of ground states is more delicate. To treat it we need the following

\begin{Theorem}\label{t4.1} In addition to assumptions (i) -- (vi), suppose that
\begin{equation}\label{e4.1}
|f(x, u+v)-f(x,u)\leq C|v|(1+|u|^{p-1})\quad\mbox{if } |v|\leq\varepsilon
\end{equation}
with some positive $C$ and $\varepsilon$.

(a) Let $u_k\in E$ be a bounded sequence such that $J'(u_k)\to 0$ and $J(u_k)\to c>0$. Then there exist critical points
$u^i$ of $J$
and sequences $b^i_k\in\mathbb{Z}^N$, $i=1,\ldots , n$ such that
\begin{equation}\label{e4.2}
\sum_{i=1}^n J(u^i)=c
\end{equation}
and
$$
\|u_k-\sum_{i=1}^n u^i(\cdot +b^i_k)\|\to 0\quad\mbox{as } k\to\infty.
$$
(b) Let $u_k\in E_k$ be a  sequence such that $\|u_k\|_k$ is bounded, $J'_k(u_k)\to 0$ and $J_k(u_k)\to c>0$. Then there exist critical points
$u^i$ of $J$ satisfying (\ref{e4.2})
and sequences $b^i_k\in\mathbb{Z}^N$, $i=1,\ldots , n$ such that
$$
\|u_k-\sum_{i=1}^n u^i(\cdot +b^i_k)\|_{H^1(Q_k)}\to 0 \quad\mbox{as } k\to\infty.
$$
\end{Theorem}

For the proof of part (a) we refer to \cite{k-sz}. Part (b) is borrowed from \cite{pank:nls}.

{}From Theorem~\ref{t4.1} and previous remarks it follows

\begin{Theorem}\label{t4.2}
Under assumptions (i) -- (vi) and (\ref{e4.1}) there exist at least one $k$-periodic ground state and at least one ground state.
\end{Theorem}

Now we want to study the behavior of $k$-periodic ground states as $k\to \infty$. We replace assumption $(v)$ by the following stronger
assumption

$(v')${\em The function $f(x, u)$ is $C^1$ in $u$, $f'_u$ is a Carath\'eodory function and there exists $\theta \in (0,1)$ such that
$$0<u^{-1}f(x,u)\leq\theta f'_u(x,u)$$
for every $u\neq 0$ and
\begin{equation}\label{e4.a}
|f'_u(x,u)|\leq C(1+|u|)^{p-2}
\end{equation}
for all $u\in\mathbb{R}$, with $C>0$.}

This implies, in particular, that
\begin{equation}\label{e4.3}
0<F(x, u)\leq\frac{\theta}{1+\theta}f(x,u)u
\end{equation}
for $u\neq 0$, and $\theta/(1+\theta)<1/2$. Therefore, $(v)$ follows from $(v')$. Moreover, in this case inequality (\ref{e4.1}) holds as well.
The functionals $J_k$ and $J$ are $C^2$, and
\begin{align}\label{e4.3a}
\langle J_k''(u)w, v\rangle &= (P_k^+w,v)-(P_k^-w,v)\pm\int_{Q_k}f'(x,u)wv\,dx=\notag\\
      &=\int_{Q_k}\left(\nabla w\nabla v +V(x)wv \pm f'(x,u)wv\right)\, dx,
\end{align}
\begin{align}\label{e4.3b}
\langle J''(u)w, v\rangle &= (P^+w,v)-(P^-w,v)\pm\int_{\mathbb{R}^N}f'(x,u)wv\,dx=\notag\\
      &=\int_{\mathbb{R}^N}\left(\nabla w\nabla v +V(x)wv \pm f'(x,u)wv\right)\, dx.
\end{align}

Define the set $S^{\pm}_k\subset E_k$ that consists of all nonzero $u\in E_k$  such that
$$I_k(u):=\langle J'_k(u),u\rangle =0$$
and
$$\langle J'_k(u),v\rangle =0\quad\forall v\in E_k^{\mp}\,.$$
Similarly, the set $S^{\pm}\subset E$ consists of all nonzero $u\in E$  such that
$$I(u):=\langle J'_k(u),u\rangle =0$$
and
$$\langle J'(u),v\rangle =0\quad\forall v\in E^{\mp}\,.$$
These sets are nonempty, because they contain solutions.

Now let us consider the following minimization problems
\begin{equation}\label{e4.4}
m^{\pm}_k= \inf\{\pm J_k(u):u\in S^{\pm}_k\}
\end{equation}
and
\begin{equation}\label{e4.5}
m^{\pm}= \inf\{\pm J(u):u\in S^{\pm}\}\,.
\end{equation}

Inspection of the proof of Lemmas~\ref{l2.1} and \ref{l2.2} gives us the following result.

\begin{lemma}\label{l4.1}
There exist constants $\varepsilon_1 >0$ and $\varepsilon_2 >0$ independent on $k$ such that
for every  $u_k\in S^{\pm}_k$ and every $u\in S^{\pm}$  we have
$\|u_k\|_k \geq \varepsilon_1$, $\|u\| \geq \varepsilon_1$, $\pm J_k(u_k)\geq\varepsilon_2$
and $\pm J(u)\geq\varepsilon_2$.Moreover,
$$\|u_k\|_k\leq C(|J_k(u_k)|^{1/2}+ |J_k(u_k)|^{1/p'})\,,\quad \|u\|\leq C(|J(u)|^{1/2}+ |J(u)|^{1/p'})\,.$$
\end{lemma}

Let $\overline{E}_k^{\pm}=\mathbb{R}\oplus E_k^{\pm}$ and $\overline{E}^{\pm}=\mathbb{R}\oplus E^{\pm}$. Consider
the operators $G_k^{\pm}:E\to \overline{E}_k^{\mp}$ and $G^{\pm}:E\to \overline{E}^{\mp}$ defined by
$$
G_k^{\pm}(u)=(\langle J'_k(u),u\rangle, P_k^{\mp}J'_k(u))\,,
$$
$$
G^{\pm}(u)=(\langle J'(u),u\rangle, P^{\mp}J'(u))\,,
$$
respectively. Obviously, $S_k^{\pm}=(G_k^{\pm})^{-1}(0)\setminus \{0\}$ and $S^{\pm}=(G^{\pm})^{-1}(0)\setminus \{0\}$. Moreover,
these operators are $C^1$, and
\begin{equation}\label{e4.5a}
(G_k^{\pm})'(u)v=\left(\langle J_k''(u)v,u\rangle +\langle J_k'(u),v\rangle\,,P_k^{\mp}J_k''(u)v\right)\,,
\end{equation}
\begin{equation}\label{e4.5b}
(G^{\pm})'(u)v=\left(\langle J''(u)v,u\rangle +\langle J'(u),v\rangle\,,P^{\mp}J_k''(u)v\right)\,.
\end{equation}

\begin{lemma}\label{l4.2}
The set $S^{\pm}_k$ (respectively, $S^{\pm}$) is a closed $C^1$-submanifold of $E_k$ (respectively, $E$).
\end{lemma}

{\em Proof\/}. We consider $S^+$ only. The remaining cases are similar. The result follows from the implicit function
theorem if we check that $G'$ is onto at every point of $S^+$. For notational convenience we skip the superscript in
$G$ and the domain of integration.

Let $u_0\in E$ ($u_0\neq 0$), $\langle J'(u_0),u_0\rangle = \alpha_0$ and $P^-\nabla J(u_0)=\alpha$. Identifying $\overline{E}^-$ and the subspace
$\mathbb{R}u_0\oplus E^-$ via $(\tau, h)=\tau u_0+h$, a direct calculation gives
\begin{equation*}\begin{split}
(G'(u_0)(\tau,h),(\tau,h))&=2\alpha_0\tau^2 -3\tau (P^-u_0,h) +\tau^2 \int [f(x,u_0)u_0-f'_u(x,u_0)u_0^2]\,dx-\\
&-(h,h)-\tau\int [2f'_u(x,u_0)u_0+f(x,u_0)]h\,dx -\int f'_u(x,u_0)h^2\, dx\,.
\end{split}\end{equation*}
Since for $h\in E^-$
\begin{equation*}\begin{split}
(P^+u_0,h) -(P^-u_0,h) -\int f(x,u_0)h\,dx&=\\
=-(P^-u_0,h) -\int f(x,u_0)h\,dx&=(\alpha,h)\,,
\end{split}
\end{equation*}
we obtain
$$
(G'(u_0)(\tau,h),(\tau,h))=2\alpha_0\tau^2 -(h,h) +3\tau (\alpha,h) -
$$
$$
-\int \{[f'_u(x,u_0)u_0^2 -f(x,u_0)u_0]\tau^2 +2[f'_u(x,u_0)u_0-f(x,u_0)]\tau h+ f'_u(x,u_0)h^2\}\,dx=
$$
$$
=2\alpha_0\tau^2 -(h,h) +3\tau (\alpha,h) -\int (A\tau^2 +2C\tau h +Bh^2)\,dx\,.
$$
Note that $B\geq 0$, and $B(x)=0$ iff $u_0(x)=0$. Hence, in the last case also $A(x)=C(x)=0$. Since
$$
A\tau^2 +2C\tau h + Bh^2=\left(A-\frac{C^2}{B}\right)\tau^2 +\left(\sqrt{B}h+\frac{C\tau}{\sqrt{B}}\right)\,,
$$
$$
|\tau (\alpha , h) |\leq\frac{1}{2}\|\alpha\|(\tau^2+\|h\|^2)\,
$$
and, due to assumption $(v')$,
\begin{equation*}\begin{split}
A-\frac{C^2}{B}&=\left(u_0-\frac{f(x,u_0)}{f'_u(x,u_0)}\right)f(x,u_0)\geq\\
&\geq (1-\theta )f(x,u_0)u_0\,,
\end{split}\end{equation*}
we obtain
\begin{equation}\label{e4.5c}\begin{split}
(G'(u_0)(\tau,h),(\tau,h))\leq & 2\alpha_0\tau^2 -\|h\|^2+\frac{3}{2}\|\alpha\|\tau^2+\frac{3}{2}\|\alpha\|\|h\|^2-\\
& -\tau^2(1-\theta)\int f(x,u_0)u_0\,dx\,.
\end{split}\end{equation}

Now if $u_0\in S^+$, then $\alpha_0=0$ and $\alpha =0$. Equation (\ref{e4.5a}) implies that on the subspace
$\mathbb{R}u_0\oplus E^-\subset E$ the operator $G'(u_0)$ is strictly negative defined, hence, invertible. Therefore,
the implicit function theorem applies and we conclude. \eproof

\begin{rem}\label{r4.1}{\em It is clear that $S^{+}_k$ has a finite co-dimension, while $S^{-}_k$ is
finite dimensional. The tangent space to $S^{\pm}_k$ at $u_0$ consists of all $h\in E_k$ such that
$$
\langle I'_k(u_0), h\rangle =0\,,\quad P^{\mp}_k \nabla J_k(u_0)h=0
$$
and similarly for $S^{\pm}$. The subspace $\mathbb{R}u_0\oplus E_k^{\mp}\subset E_k$ (resp., $\mathbb{R}u_0\oplus E^{\mp}\subset E$) is
transverse to $S^{\pm}_k$ (resp., to $S^{\pm}_k$) at $u_0$, as it follows from the proof of Lemma~\ref{l4.2}.}
\end{rem}

\begin{lemma}\label{l4.3}
Any critical point of the restriction of $J_k$ to $S_k$ (respectively, $J$ to $S$) is a critical point of $J_k$
(respectively, $J$). In particular, solutions of problems (\ref{e4.4}) and (\ref{e4.5}) are periodic ground states and
ground sates, respectively.
\end{lemma}

{\em Proof\/}. A direct calculation shows that if $u_0\in S_k^{\pm}$ (resp., $u_0\in S^{\pm}$), then $J'_k(u_0)$ (resp., $J'_k(u_0)$)
vanishes on the subspace $\mathbb{R}u_0\oplus E_k^{\mp}\subset E_k$ (resp., $\mathbb{R}u_0\oplus E^{\mp}\subset E$). Therefore, if
$u_0$ is a critical point of $J_k$ (resp., $J$) restricted to $S_k^{\pm}$ (resp., to $S^{\pm}$), then $J'_k(u_0)$ (resp., $J'_k(u_0)$)
vanishes everywhere, since it is equal to 0 also on the tangent space. \eproof

\begin{Theorem}\label{t4.3}Problem~(\ref{e4.4}) has at least one solution which is a periodic ground state. Moreover, $m_k^{\pm}\geq
\varepsilon_2$, with $\varepsilon_2$ from Lemma~\ref{l4.1}.
\end{Theorem}

{\em Proof\/}. We only sketch the existence of the minimization problem on $S_k^{+}$. On this set we have
\begin{equation}\label{e4.5d}
 J_k(u)= J_k(u)-\frac{1}{2}\langle J'_k(u),u\rangle=\int_{Q_k}\left(\frac{1}{2}f(x,u)u-F(x,u)\right)\,dx\,.
\end{equation}
Hence, due to assumption $(v)$, $J_k(u)$ is bounded below on $S_k^{+}$. Consider a minimization sequence $u_n\in S_k^{+}$. As it follows from the
Ekeland principle (see, e.g., \cite{zeidler}), we can assume that $u_n$ is, in addition, a Palais-Smale sequence for ${J_k}_{|S^+_k}$.
(In \cite{zeidler} functionals defined on a Banach space are considered, but arguments work in the case of functionals on $C^1$ Banach manifolds,
with only minor change). We shall show that, in fact, $u_n$ is a Palais-Smale sequence for the whole functional $J_k$. This is enough to pass to the
limit and get the result.

Let $g_n=\nabla J_k(u_n)$ and $g^t_n$ its component tangent to $S^+_k$, i.e. the orthogonal projection of $g_n$ onto the tangent space at $u_n$.
Then $g^t_n\to 0$. We have to show than $g_n\to 0$. Inequality (\ref{e4.5c}), with $u_0=u_n$, $\alpha_0=0$ and $\alpha=0$, implies that the operator
${(G^+_k)}'(u_n)$ has a right inverse operator, $A^n_k$. The image of this operator is the subspace $\overline{E}^-(u_n)=\mathbb{R}u_n\oplus E^-$ and
its norm is bounded above by a constant $C>0$ independently of $n$. The operator
$A^n_k{(G^+_k)}'(u_n)$ is the projector, $P^n_k$, onto $\overline{E}^-(u_n)$ parallel to
the tangent space at $u_n$. The adjoint operator $(P^n_k)^*$ is the projector onto the orthogonal complement to $\overline{E}^-(u_n)$ parallel to
the normal subspace at $u_n$. Hence, $g_n=(P^n_k)^*g_n^t$ and $\|g_n\|_k\leq C\|g_n^t\|$. This completes the proof. \eproof

Now we impose the following additional assumption

$(vii)${\em There exist $C>0$ and $\gamma\in (0,1]$ such that
$$
|f'_u(x,u_1)-f'_u(x,u_2)|\leq C(1+|u_1|+|u_2|)^{p-2-\gamma}|u_1-u_2|^{\gamma}
$$
for all $u_1,u_2\in\mathbb{R}$.}

This assumption implies, in particular, (\ref{e4.a}).

\begin{Theorem}\label{t4.4}Under assumptions $(i)$--$(iv)$, $(v')$, $(vi)$ and $(vii)$ we have
\begin{equation}\label{e4.6a}
m^{\pm}=\lim_{k\to\infty}m^{\pm}_k\,.
\end{equation}
Moreover, let $u_k\in E_k$ be a solution of (\ref{e4.4}). Then, after passage to a subsequence still denoted by $u_k$,
there exist a solution $u$ of (\ref{e4.5}) and a sequence $b_k\in\mathbb{Z}^N$ such that
\begin{equation}\label{e4.6}
\|u_k-u(\cdot +b_k)\|_{H^1(Q_k)}\to 0\,.
\end{equation}
\end{Theorem}

{\em Proof\/}. We again consider only the case of ``+'' sign. Since $u_1$ is also a $k$-periodic solution, we see that $m^+_k\leq m^+_1$. Hence,
the sequence $m^+_k$ is bounded and Lemma~\ref{l4.1} shows that the sequence $\|u_k\|_k$ is bounded. Due to Theorem~\ref{t3.1} we can assume that
$u_k$ converges in $H^1_{\mbox{loc}}(\mathbb{R}^N)$ to a nontrivial solution $u\in E$. Note that $J(u)\geq m^+$.

Equation (\ref{e4.5d}) implies that
$$
m^+_k=\int_{Q_k}\left(\frac{1}{2}f(x,u_k)u_k-F(x,u_k)\right)\,dx=:\int_{Q_k}g(x,u_k)\,dx\,.
$$
Since the integrand here is nonnegative, we see that for any bounded domain $\Omega\subset\mathbb{R}^N$ and $k$ large enough
$$
m^+_k\geq\int_{\Omega}g(x,u_k)\,dx
$$
and, therefore,
$$
\liminf m^+_k\geq\int_{\Omega}g(x,u_k)\,dx\,.
$$
Since $\Omega$ is an arbitrary domain and $u$ is a nontrivial solution, we obtain
\begin{equation}\label{e4.7}
\liminf m^+_k\geq\int_{\mathbb{R}^N}g(x,u_k)\,dx\,=J(u)\geq m^+\,.
\end{equation}

To prove (\ref{e4.6a}), we have now to show that
\begin{equation}\label{e4.7a}
\limsup m^+_k\leq m^+\,.
\end{equation}
Let $v\in S^+$. Since $C^{\infty}_0(\mathbb{R}^N)$ is dense in $E$ there exists a sequence of smooth functions $v_k\to
v$ in $E$ such that $v_k$ has a compact support in $Q_k$. Having in mind periodic extension, we can also consider
$v_k$ as an element of $E_k$. Let $\alpha^k=P^-_k\nabla J_k(v_k)$ and $\alpha^k_0=I_k(v_k)$. Inequality (\ref{e4.5c}),
with $u_0=v_k$, shows that the derivative of the map $G^+_k$ at $v_k$ has a right inverse whose norm is bounded
independently of $k$, provided $k$ is large enough. Note that, due to assumption $(vii)$, $G^+_k$  is H\"older
equicontinuous, with the exponent $\gamma$, say, on the ball of radius $2\|v\|$ in $E_k$ centered at the origin.
Inspecting the proof of the implicit function theorem (see, e.g., \cite{zeidler}), we get right inverse map, $T_k$, to
$G^+_k$ defined on the ball $B_k\subset \mathbb{R}\oplus E^-_k$ of radius $r$ centered at $(\alpha_0^k,\alpha^k)$,
where $r$ is independent of $k$. Moreover, $T_k$ is H\"older eqicontinuous. For $k$ large enough, $0\in B_k$ and is
close to the point $(\alpha_0^k,\alpha^k)$. Hence, $w_k=T_k(0)$ is close to $v_k$, $w_k\in S^+_k$ and
$$
J_k(w_k)\leq J(v)+\varepsilon
$$
if $k$ is large enough. This implies (\ref{e4.7a}).

The remaining part of the theorem follows from Theorem~\ref{t4.1}, part (b). \eproof

\section{Exponential decay}\label{s5}

Now let us study the decay of finite action solutions to equation (\ref{e1}).

\begin{Theorem}\label{t5.1}
Let $u\in H^1(\mathbb{R}^N)$ be a solution of (\ref{e1}). Then $u$ is a continuous function and
there exist positive constants
$C$ and $\lambda$ such that
$$
|u(x)|\leq C\exp(-\lambda |x|)\,
$$
where $\lambda \leq c\,\mbox{dist}\,(0,\sigma (L))= c\alpha$.
\end{Theorem}

{\em Proof\/}. First, we have to show that $u\in L^{\infty}(\mathbb{R}^N)$. This can be done
exactly as in the proof of Lemma~5.1, \cite{b-p-w}. (That proof is based on the Sobolev estimates
for Schr\"odinger operators). Next, Theorem~B.3.3 of \cite{si} implies that $u$ is continuous.
Now we set $W(x):=-f(x,u(x))/u(x)$ (with
$W(x)=0$ if $u(x)=0$). Hence, $u$ solves the equation
\begin{equation}\label{e5.1}
-\Delta +(V(x)+W(x))u=0
\end{equation}
on $\mathbb{R}^N$. Since $V+W\in L^{\infty}(\mathbb{R}^N)$, Theorem~C.3.1 of \cite{si} shows that
$u(x)\to 0$ as $x\to \infty$.

Since $W(x)\to 0$ as $x\to\infty$, the potential $W$ defines a relatively compact perturbation of the operator $L$.
Hence, $L+W$ has the same essential spectrum as $L$ and may have only isolated eigenvalues of finite multiplicity
outside the essential spectrum. Equation (\ref{e5.1}) means that 0 is an eigenvalue  of $L+W$ outside the essential
spectrum, hence, of finite multiplicity, and $u$ is the corresponding eigenfunction. Now Theorem~C.3.5, \cite{si},
provides the required exponential bound.\eproof

\begin{rem}\label{r5.1}{\em
Using estimates for Green's function of $L$ (see \cite{kuch}, Section~7.6.3), one can improve the estimate for $\lambda$ as follows:
$a\leq c\alpha^{1/2}$.}
\end{rem}

\section{Further results}\label{s6}

Now we discuss assumption $(vi)$. More precisely, what may happen if $0$ {\em is} in the spectrum of $L$?
It is commonly believed that in this case well localized solutions of (\ref{e1}) do not exist.
As we see from equation (\ref{e5.1}), the question
is closely related to the problem of absence of embedded eigenvalues for periodic Schr\"odinger operators
perturbed by a decaying potential. Unfortunately, there is no satisfactory general result on embedded
eigenvalues, except the case $N=1$. In the case $N=1$ F.~S.~Rofe-Beketov \cite{r-b} (see also \cite{r-b-k})
has proved that $L+W$ has no embedded eigenvalues if $(1+|x|)W(x)\in L^1(\mathbb{R})$. This
implies immediately the following result.

\begin{Theorem}\label{t6.1}
Let $N=1$. Suppose that assumptions $(i)$ and $(ii)$ are satisfied and, for some $r>1$,
$|f(x,u)|\leq c |u|^r$ for a.e. $x\in\mathbb{R}$ and $u$ in a neighborhood of $0$. Then equation
(\ref{e1}) has no nontrivial weak solutions such that $(1+|x|)|u(x)|^{r-1}\in L^1(\mathbb{R})$. In
particular, there is no nontrivial exponentially decaying solutions.
\end{Theorem}

On the other hand, T.~Bartsch and Y.~Ding have found the following existence result.

\begin{Theorem}\label{t6.2}
In addition to $(i)$--$(v)$, assume that $V$ and $g$ are continuous functions of all their arguments,
there exist  constants $r\in [q, 2^*)$, $s\in [p, 2^*)$,  $a_1>0$ and $a_2>0$ such that
$$a_1|u|^r\leq F(x,u)\,,$$
$$|f(x,u)|\leq a_2(|u|^{p-1}+|u|^{s-1})$$
for all $x\in\mathbb{R}^N$ and $u\in\mathbb{R}$. Let $0\in\sigma(L)$ and there exists $b >0$ such
that $(0,b)\cap\sigma(L)= \emptyset$ in the case of `$+$' sign in (\ref{e1})
and $(-b, 0)\cap\sigma(L)= \emptyset$
in the case of `$-$' sign. Then equation (\ref{e1}) has a nontrivial weak solution
$u\in H_{\mathrm{loc}}^2(\mathbb{R}^N)$. This solution lies in $L^t(\mathbb{R}^N)$ for
$r\leq t\leq 2^*$, is continuous and $u(x)\to 0$ as $x\to\infty$.
\end{Theorem}

Certainly, the solution obtained in the last theorem should be not well decaying. If $N=1$, Theorem~\ref{t6.1}
implies that this solution cannot decay too fast.

Thus, we conjecture that if $0\in\sigma (L)$, then nontrivial well-decaying at infinity solutions to periodic NLS
\ref{e1} do not exist. To support this conjecture assume that $0\in\sigma (L)$ and is not a lower (resp., upper) edge
of a spectral gap in the case of `+' (resp., `-') sign in (\ref{e1}). Let $k_j\to infty$ be a sequence of integers
such that $0\not\in\sigma (L_{k_j})$. Recall that $\sigma (L_k)$ is discrete and the countable set $\cup\sigma (L_k)$
is dense in $\sigma (L)$. Therefore, $k_j$ is a generic sequence. Let $u_j=u_{k_j}\in E_{k_j}$ be the $k_j$-periodic
solution obtained by linking (see Remark~\ref{r2.*}) and $c_j=c_{k_j}$ be the corresponding linking critical value
defined in Theorem ~\ref{t2.1}. Note that it does not depend on particular choice of $z^0$. Denote by $\alpha_j$ the
distance from 0 to $\sigma (L_{k_j})$. Then we have
\begin{prop}\label{p6.1}
Under the assumptions imposed above
$$
\|u_j\|_{H^1(Q_{k_j})}\leq C\alpha^{\frac{q}{2(q-2)}}\,,
$$
where $C>0$ does not depend on $j$.
\end{prop}

{\em Proof\/}. Since the result seems to be not exact, we only sketch the proof in the case of `+' sign.

Following the proof of Lemma~\ref{l2.1} (see also Remark~\ref{r2.a}), we obtain
\begin{equation}\label{e6.*}
\|u_j\|_{H^1(Q_{k_j})}\leq C\alpha^{-1/2}(|c_j|^{1/p'}+|c_j|^{1/2})\,,
\end{equation}
where $p'$ is the conjugate exponent to $p$.

Now we estimate the linking critical  value as in the proof of Theorem~\ref{t2.1}. We refine the argument after inequality (\ref{e2.5}).
Take as $z^0\in E^+_{k_j}$, the eigenvector of $L_{k_j}$ with $\|z^0\|_{k_j}=1$ that correspond to the first positive eigenvalue $\lambda_0$
of $L_{k_j}$. Then $\|z^0\|_{H^{-1}(Q_{k_j})}\leq C\|z^0\|_{L^q(Q_{k_j})}$, with $C>0$ independent of $j$. Since
$$
-\Delta z^0 + z^0= \lambda_0 z^0 -V(x)z^0\,,
$$
we have that $\|z^0\|_{H^1(Q_{k_j})}\leq C\|z^0\|_{L^q(Q_{k_j})}$. From (\ref{e2.5}) we get
$$
c_j\leq\sup_MJ_{k_j}(u)\leq \max (3t^2-C\alpha_j^{-q/2}t^p)\,.
$$
A direct calculation gives us
$$
c_j\leq C\alpha_j^{\frac{q}{q-2}}\,.
$$
Since $1<p'<2$, this together with (\ref{e6.*})  implies the required. \eproof

\begin{rem}{\em Certainly, the conclusion of Proposition~\ref{p6.1} holds also for periodic ground states.}
\end{rem}

Construction of ground states in Section~{s4} shows that the ground critical value does not exceed the linking value, which is not necessary
a critical value in this case. Therefore, arguing as in Proposition~\ref{p6.1}, we obtain the following result.

\begin{prop}\label{p6.2} Let $u\in E$ be a ground state of equation (\ref{e1}). Under the assumptions of Theorem~\ref{t4.2} we have the estimate
$$
\|u\|_{H^1(\mathbb{R}^N)}\leq C(\alpha_{\pm})^{\frac{q}{2(q-2)}}\,,
$$
where $(-\alpha_-,\alpha_+)$ is the spectral gap containing 0.
\end{prop}

As consequence, we see that ground states bifurcate from the trivial solution corresponding to an appropriate edge of the spectral gap.

\section{Gap solitons}\label{gs}

In this section we give an application of the previous results to a problem that arises in the
theory of photonic crystals, the existence of gap solitons. Photonic crystals are dielectric media
with spatially periodic (or close-to-periodic) structure.
A good introduction into this field can be found in \cite{j-m-w}. For a survey of rigorous mathematical
results we refer to \cite{kuch}. Both these publications deal with linear optical media. If the medium we
consider is nonlinear, many new phenomena occur. Among them one of the most interesting is the possibility
of gap solitons, i.e. spatially localized light patterns with the frequency prohibited by the linear theory.
On physical level a simple description of this phenomenon is presented in \cite{mills}. However, up to now
there was no mathematically rigorous proof of the existence of gap solitons even in simplest situations. Here
we show that some existence result can be extracted from the results on periodic NLS.

In a dielectric medium we consider the system of Maxwell equations
$$
\nabla\times\mathbf{E}=-\frac{\partial\mathbf{B}}{\partial t}\,,
$$
$$
\nabla\times\mathbf{H}=\frac{\partial\mathbf{D}}{\partial t}\,,
$$
$$
\nabla \mathbf{D}=0\,,
$$
$$
\nabla\mathbf{B}=0\,.
$$
Since we assume that the medium is non-magnetic, we set $\mathbf{B}=\mathbf{H}$. We consider the following
constitutive relation between displacement and electric fields:
$$
\mathbf{D}=(\varepsilon (x) +\chi (x)\langle|\mathbf{E}|\rangle^2)\mathbf{E}\,
$$
where $\langle\cdot\rangle$ stands for time average. Such type of nonlinear response was introduced by
N.~N.~Akhmediev \cite{akh} (see also \cite{stu} and references therein). The cubic form of the nonlinearity
means that we concentrate on Kerr-like media.

In the case of two dimensional structure we assume that the functions $\varepsilon (x)=\varepsilon (x_1,x_2)$
and $\chi (x)=\chi (x_1,x_2)$ are independent on $x_3$, periodic in $(x_1,x_2)$ and of the class
$L^{\infty}$. Moreover, we suppose that $\varepsilon (x)\geq\varepsilon_0 > 0$ and either
$\chi (x) >0$, or $\chi (x)<0$ everywhere. The assumption means that the medium we consider is either
everywhere self-focusing, or everywhere defocusing \cite{mills}. We restrict ourself to the so-called
$E$-mode $\mathbf{E}=(0,0,E)$ and look for solutions of the form
\begin{equation}\label{e-gs1}
E=u(x_1,x_2)\cos (\beta x_3-\omega t +\theta_0)
\end{equation}
where the amplitude vanishes at infinity. Such solutions represent  light patterns of the frequency $\omega$
that are localized in $(x_1,x_2)$ directions and propagate along the $x_3$-axis, with the wave number $\beta$.
Obviously, the problem reduces to the two dimensional NLS
\begin{equation}\label{e-gs2}
-\Delta u-\omega^2\varepsilon u+\beta^2u=\omega^2\chi u^3.
\end{equation}
Now we can apply the previous results. In the self-focusing case ($\chi >0$) we obtain a nontrivial
exponentially decaying solution if $-\beta^2$ is not in the spectrum of the operator
$L_{\omega}:=-\Delta-\omega^2\varepsilon$. In particular, given $\omega$ such a solution exists
for all $|\beta|$ large enough. In the defocusing case ($\chi <0$) nontrivial exponentially localized solutions
exist, provided $-\beta^2$ is not in $\sigma (L_{\omega})$, but {\em not} below the spectrum. This means that
$-\beta^2$ must belong to a {\em finite} spectral gap of $L_{\omega}$.

The case  $\beta =0$ is of particular interest. Let us assume that $0\not\in\sigma (L_{\omega})$. This means
that $\omega$ is a prohibited frequency \cite{kuch}. Since the potential of $L_{\omega}$ is {\em negative},
$0$ belongs to a finite gap. Therefore, independently of the sign of $\chi$ there exists a nontrivial solution
of the form $u(x_1,x_2)\cos (-\omega t +\theta_0)$ with exponentially decaying amplitude $u$. Such solutions
represent so-called {\em standing} gap solitons. Hence, we obtain the existence of gap solitons for every
prohibited frequency. In fact, in this case ($0\in\sigma(L_{\omega})$) there exists a nontrivial exponentially
localized solution of the form (\ref{e-gs1}) for all $\beta$ close enough to $0$.

Exactly the same results hold for one dimensional structures. This case reduces to the one dimensional
periodic NLS.

Now let us discuss the behavior of gap solitons with respect to $\omega$. Let $(-\alpha_-, \alpha_+)$ be the spectral
gap of $L_{\omega}$ containing 0 and $(\omega_-,\omega_+)$ the corresponding gap of (nonnegative) frequencies. Note
that if $\omega$ goes to $\omega_{\pm}$, then $\alpha_{pm}$ goes to 0. Indeed, the spectrum of $L_{\omega}$ is the
union of closed intervals
$$[\min_{\theta\in [0,2\pi]}\lambda_j(\theta),\max_{\theta\in [0,2\pi]}\lambda_j(\theta)]\,,$$
where $\lambda_j(\theta)$ are the Bloch eigenvalues (see, e.g., \cite{r-s}). Since the potential $-\omega^2\varepsilon (x)$ decreases monotonically
as $\omega$ increases, the comparison principle for eigenvalues (see, e.g., \cite{b-sh}) tells us that the eigenvalues $\lambda_j(\theta)$ decreases
monotonically and continuously, and we are done. Now Proposition~\ref{p6.2} shows that gap solitons bifurcate from zero solutions corresponding to
$\omega=\omega_+$ in the  self-focusing case and $\omega=\omega_-$ in the defocusing case.

We have considered here the case when $\chi (x)$ does not change sign. The case of sign changing $\chi (x)$ is not less important, but
completely open. It corresponds to a mixture of self-focusing and defocusing optical materials.

Another important problem is the existence of gap solitons in periodic media with {\em saturation\/}. In this case the nonlinearity
is asymptotically linear. A particular examples of such nonlinearities are
$$
f(x,u)=\chi (x)u^2(1 + c(x)u^2)^{-1}\,
$$
$$
f(x,u)=\chi (x)(1-\exp (-a(x)u^2))\,.
$$
In the context of applications to nonlinear optical wave guides (one dimensional nonperiodic problem) such
nonlinearities were studied extensively by C.~Stuart (see, in particulary, \cite{stu}).

\end{document}